\documentclass[12pt,a4paper]{article}
\usepackage{geometry}
\usepackage{graphicx}
\usepackage{amsfonts,amssymb,eucal}
\usepackage{booktabs}

\usepackage[reqno]{amsmath}

\newtheorem{TT}{Theorem}
\newtheorem{CC}{Corollary}
\newtheorem{Conj}{Conjecture}
\newtheorem{RR}{Remark}

\begin{document}

\title{\bf Number of Vertices of the Polytope \\
of Integer Partitions and Factorization \\ of the Partitioned Number%
\footnote{Dedicated to the memory of Professor V. B. Priezzhev}}

\author{\Large Vladimir A. Shlyk%
\footnote{Belarusian State University, v.shlyk@gmail.com}}

\date{\vspace{-5ex}}

\maketitle

\begin{abstract}
The polytope of integer partitions of $n$ is the convex hull of the corresponding $n$-dimensional integer points.
Computation shows intriguing features of $v(n),$ the number of the polytope vertices:
its graph has a tooth-shaped form with the highest peaks at prime $n$'s.
We explain the shape of $v(n)$ by the large number of partitions of even $n$'s that were
counted by N. Metropolis and P. R. Stein in 1970.
We reveal that divisibility of $n$ by 3 also reduces the value of $v(n),$
which is caused by partitions that are convex combinations of three but not two others,
and characterize convex representations of such integer points in arbitrary integral polytope.
We use a specific classification of integers and
demonstrate that the graph of $v(n)$ is stratified into layers corresponding to resulting classes.
Our main conjecture claims that $v(n)$ depends on small divisors of $n.$
We also offer an argument for that the numbers of vertices of the master corner polyhedron on the cyclic group
have features similar to those of $v(n).$
\end{abstract}

\section{Introduction}

Integer partitions are related to divergent problems in mathematics and statistical mechanics \cite{And76}.
A partition of a positive integer $n$
is any finite non-decreasing sequence $\rho $ of positive integers $n_1, n_2, \ldots, n_r$ such that
\begin{equation*}
 \sum_{j=1}^r n_j = n.
\end{equation*}
The integers $n_1, n_2, \ldots, n_r$ are called parts of the partition $\rho.$

In this paper we develop the polyhedral approach to integer partitions proposed in \cite{Shl1}.
It is based on the $n$-dimensional geometrical interpretation of partitions \cite{WeisWolfram}.
Every partition $\rho$ is referred to as a non-negative integer point
$x=(x_1,x_2,\ldots , x_n)\in \mathbb{R}^n,$ a solution to the equation
\begin{equation}\label{eqPn}
 x_1 + 2x_2+ \ldots +nx_n =n,
\end{equation}
with $x_i,$ $i=1, \ldots, n,$ being the number of parts $i$ in $\rho.$
For example, the partition $8=4+2+1+1$ with three distinct parts $1,$ $2,$ $4$ is identified with
$x=(2,1,0,1,0,0,0,0)\in\mathbb{R}^8.$
We keep on writing $x\vdash n$ to indicate that $x\in\mathbb{R}^n$ is a partition of $n.$

Let $P(n)$ denote the set of partitions of $n.$
The polytope of partitions of $n,$ $P_n\subset \mathbb{R}^n,$ is defined as the convex hull of $P(n):$
$$
P_n := \mathrm{conv\,}P(n)= \mathrm{conv\,}\{x=(x_1,x_2,\ldots,x_n)\in \mathbb{R}^n ~|~x\vdash n\}.
$$
The conversion from the set to a polytope reveals the geometrical structure of $P(n).$
As for every polytope, the key elements of $P_n$ are its facets and vertices.
The facets were characterized in \cite{Shl1}, the vertices were studied in \cite{Shl08, Shl13, Shl2}.
Vertices of $P_n$ and their number are of special importance since, by Carath\'eodory's theorem \cite{Cara},
every partition is a convex combination of some vertices.
We computed vertices of $P_n$ for $n\le 100,$ see \cite{Vroubl}, and presented their numbers
in the On-Line Encyclopedia of Integer Sequences (OEIS), sequence A203898.
It turned out that the number of vertices of $P(n)$ is much less than the number of partitions of $n.$

The problem of recognizing vertices of $P(n)$ is proved to be decidable in polynomial time with the use of
linear programming technique \cite{Onn}.
However, no combinatorial characterization of vertices of $P_n$ is available as yet.
The only result in this direction
is the criterion for a partition to be a convex combination of two others \cite{Shl13}, see Theorem \ref{Criterion} below.

Let ${\rm Vert}\, P_{n}$ denote the set of vertices of $P_n$ and $v(n):={\rm |Vert}\, P_{n}|$ be the number of vertices.
The graph of $v(n)$ exhibits peculiar features, see Fig. 1.
In contrast to $p(n)=|P(n)|,$ the number of partitions of $n,$
the function $v(n)$ does not increase monotonously.
It drops down at every even $n$ and its peaks at prime $n$'s seem to be higher than others.
Inspired by these perplexing peculiarities, we concentrate on the
asymptotic dependence of $v(n)$ on the multiplicative properties of $n.$

\begin{center}
\noindent \includegraphics*[width=6.3in, height=3.0in,
keepaspectratio=false]{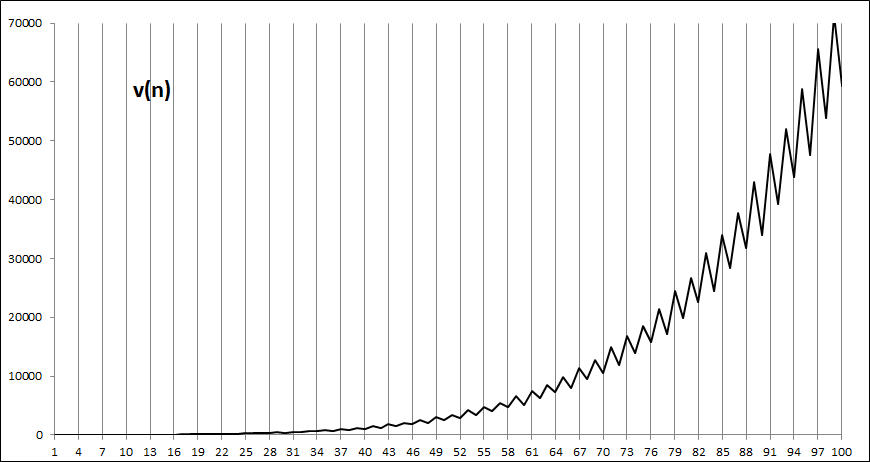}

\vskip2ex
\noindent {\bf Fig. 1.} Graph of the function $v(n),$ number of vertices of $P_{n}.$

\end{center}

To analyse the observed phenomena we study the structure of the set
of partitions of $n$ that 'fail' to be vertices of $P_n.$
In Section 2, we divide this set to subsets of partitions $x\vdash n$ according to the minimal number of partitions
that are needed to express $x$ as their convex combination.
We characterize the coefficients in such representations in the case when $x$ needs three partitions.
In fact, we do that for analogous integer points in an arbitrary integral polytope.
In Section 3, we establish that the majority of partitions in $P(n) \setminus {\rm Vert}P_{n}$
belong to a class of partitions
that were counted by N. Metropolis and P. R. Stein \cite{Metro}.
Gradually, from the obvious dependence of $v(n)$ on the evenness of $n,$
it becomes clear that other divisors of $n$ also affect the value of $v(n).$
In Section 4, we use the known numbers of knapsack partitions
introduced by R. Ehrenborg and M.A. Readdy \cite{EhRe07},
to disclose dependence of $v(n)$ on divisibility of $n$ by 3.

We suggest several conjectures grounded on our computational experience, numerical data, or some visual but quite convincing arguments.
Of primary interest is Conjecture 4 in Section 5.
In short, it claims that the value of $v(n)$ depends on factorization of $n.$
In more details, we classify integer numbers $n$ by their proximity, in a certain sense, to primes,
which is determined by the smallest divisors of $n.$
The conjecture claims that the graph of $v(n)$ has the structure of a layered cake:
it is stratified into layers corresponding to the introduced classes.
The topmost layer corresponds to the class of primes.
For a given $n,$ the major influence on the height of the layer $v(n)$ belongs to
is rendered by the smallest divisor of $n$ and every successive divisor makes its additional contribution to lowering it.
We visually demonstrate this stratification in a fragment of the graph of $v(n)$ in Figure 5.

The polytope of integer partitions has the famous 'elder brother',
the master corner polyhedron on a group.
Discovered by R. E. Gomory \cite{Gom} in 1969,
it has proved to be of the key importance in integer programming.
Gomory calls it the atom in this field \cite{Atom}.
Regrettably, after the pioneering work \cite{Gom}
its vertices fell out of research for about 40 years in pursuit of efficient cuts for integer programs
that are induced by its facets.
In Section 6, we offer an initial visual argument that the numbers of vertices of the
master corner polyhedron on the cyclic group have features similar to those of $v(n).$

Finally, in Section 7, we outline the most promising directions for the future study.

\section{Convex representations of non-vertices}

For arbitrary polytope $P,$ a point
$x\in P$ is a vertex of $P$ if it cannot be expressed as a convex combination $x=\sum_{j=1}^k \lambda_j
y^j,$ $\sum_{j=1}^k \lambda_j =1,$ $\lambda_j >0,$ of some other points $y^j \in P,$
$j=1,\ldots, k,$ in particular, of vertices.
So, every partition $x \in P(n) \setminus {\rm Vert}P_{n}$ is a convex combination of some partitions of $n.$
Denote by $\xi(x)$ the minimal number of partitions of $n,$
which are needed for such a representation of an $x{\rm \,\vdash\; }n,$ $x\notin {\rm Vert}\, P_{n},$
and let $C_{\xi }(n)$ be the set of partitions $x{\rm \,\vdash\; }n,$ for which $\xi (x)=\xi .$
It is easy to see that the sets $C_{\xi } (n),$ $\xi =2,3,4,...,$ are pairwise different and
\begin{equation} \label{GrindEQ__1_}
{\rm Vert}\, P_{n} =P(n)\backslash \bigcup _{\xi \ge 2} C_{\xi } (n).
\end{equation}

While computing vertices of $P_{n},$ we saw that for all $n$ the most of
$x{\rm \,\vdash\; }n,$ $x\notin {\rm Vert}\, P_{n} ,$
are convex combinations of some two partitions of $n,$ i. e., these $x$ belong to $C_{2} (n).$
The following theorem gives a criterion for a partition $x{\rm \,\vdash\; }n$ to belong to $C_{2} (n).$

\begin{TT}{{\rm(\cite{Shl08})}} \label{Criterion}
A partition $x{\rm \,\vdash\; }n$ is a convex combination of two partitions of $n$
if and only if there exist two different collections of parts of $x$ with equal sums.
\end{TT}

For  $n<15,$ all partitions $x\notin {\rm Vert}\, P_{n}$ belong to $C_{2} (n).$
For $n=15,21,24,25,27,28$ and $n\ge 30$ there exist non-vertices of $P_{n}$ that belong to $C_{3} (n)$
(and hence do not belong to $C_{2}(n)$).
The partition $x=(0,0,2,1,1,0^{10})\vdash 15$ corresponding to $15=3+3+4+5$ is an example;
here $0^{10}$ stands for 10 zeros.
Indeed,
$x = \frac13(0,0,5,0^{12}) + \frac13(0,0,1,3,0^{11}) + \frac13(0,0,0,0,3,0^{10}),$
there are no other partitions of 15 with parts 3, 4, 5, and $x$ is not a convex combination of any two of these partitions.

The minimal $n$ for which some $x{\rm \,\vdash\; }n$ belongs to $C_{4}(n)$ is $n=36.$
This is the partition $36=7+8+9+12,$ which is one quarter of the sum of partitions
$7^4+8, 8^3+12, 9^4, 12^3.$

So, for $n$ sufficiently large, $C_{2} (n),\, C_{3} (n),\, C_{4} (n) \neq \emptyset.$
Non-emptiness of $C_{\xi } (n)$ for $\xi \ge 5$ is not confirmed yet but we dare to suggest the following conjecture.

\begin{Conj}
For every $\xi,$ $C_{\xi } (n)\ne \emptyset$ for sufficiently large $n>n_0(\xi ).$
For every $n,$ $|C_{\xi } (n)|$ decreases when $\xi$ grows.
\end{Conj}

If this conjecture is true then the union in (\ref{GrindEQ__1_}) can consist of arbitrarily large number of sets.
The following theorem gives an upper bound for $\xi$ such that $C_{\xi}(n)\neq \emptyset.$

\begin{TT}
If for some $n$ and $\xi>2,$ $C_{\xi}(n)\neq \emptyset$ then $\xi\leq \log_2(n+1) + 1.$
\end{TT}

\noindent
\textbf{Proof.}
Let a partition $x\in C_\xi(n),$ $\xi>2,$ have $m$ parts $\{n_1, n_2, \ldots, n_m\}.$
It is proved in \cite{Shl13} that if $m>\log_2(n+1)$ then $x\in C_2(n).$
Hence $m\leq\log_2(n+1)$ since $C_2(n)\cap C_{\xi}(n) = \emptyset.$

Let $x$ be a convex combination of $y^1, y^2, \ldots, y^\xi \vdash n.$
Then $y^j_{i}=0$ for $i\notin\{n_1, n_2, \ldots, n_m\},$ $j=1, 2, \ldots, \xi.$
Since $x\notin C_k(n)$ for $k<\xi,$ then $y^1, y^2, \ldots, y^{\xi}$
are vertices of some $(\xi-1)$-dimensional simplex in $\mathbb{R}^n$ and are affinely independent.
Then the matrix with the rows $(y^j_{n_1}, y^j_{n_2}, \ldots, y^j_{n_m}),$ $j=1, 2, \ldots, \xi-1,$
is of rank $\xi-1$ and therefore $\xi-1\leq m.$
The two inequalities imply that $\xi\leq \log_2(n+1)+1.$ \hfill $\Box$

Figure 2 shows the structure of the set $P(n)$ provided Conjecture 1 holds.
The whole rectangle corresponds to all partitions of $n.$
Vertices of $P_{n} $ form the utmost right rectangle.
The inner rectangles in order from left to right correspond to
$C_{2} (n),$ $C_{3} (n),$ $C_{4} (n), \dots , C_{k} (n),$ where $k$ depends on $n.$
The set $M_2(n),$ whose definition will be given in Section 3, consists of two parts:
a subset of $C_2(n)$ depicted as the large rectangle from the left edge to the dashed line,
and a small subset of vertices forming a tiny rectangle at the bottom right of the picture.
The set $K(n)=P(n)\setminus C_2(n)$ will be considered in Section~4.

\begin{center}
\linethickness{0.25mm}
\begin{picture}(365,80) \label{fig.2}
\scriptsize
\put(0,14){\line(1,0){380}} \put(0,64){\line(1,0){380}} \put(0,14){\line(0,1){50}} \put(380,14){\line(0,1){50}}
\put(350,14){\line(0,1){50}}        
\put(282,14){\line(0,1){50}}         
\linethickness{0.18mm}
\put(282,69){\line(0,1){6}}
\put(380,69){\line(0,1){6}}
\put(282,75){\line(1,0){33}}
\put(339,75){\line(1,0){41}}
\put(317,73){$K(n)$}
\put(346,14){\line(0,1){50}}         
\put(340,14){\line(0,1){50}}
\put(316,14){\line(0,1){50}}
\put(256,14){\line(0,1){4}}
\put(256,20){\line(0,1){4}}   \put(256,26){\line(0,1){4}}    \put(256,32){\line(0,1){4}}         \put(256,38){\line(0,24){4}}
\put(256,55){\line(0,1){3}}   \put(256,60){\line(0,1){4}}
\put(365,22){\line(0,1){4}}     \put(370,22){\line(0,1){4}}   
\put(365,22){\line(1,0){5}}     \put(365,26){\line(1,0){5}}
\put(368,22){\line(-2,-3){13}}   \put(355,2){\line(1,0){5}}     \put(362,0){$M_2(n)$}
\put(70,84){\begin{picture}(40,40)

\put(281,-38){${\rm Vert}\, P_{n}$}
\put(247,-38){$C_4(n)$}
\put(217,-38){$C_3(n)$}
\put(174,-38){$C_2(n)$}
\put(6,-55){$M_2(n)$}

\end{picture}}
\end{picture}
\end{center}

\begin{center}
{\bf Fig. 2.}  Conjectured structure of the set of partitions of $n.$
\end{center}

\begin{RR} $\xi(x)$ can be defined for any integer point $x$ in any integral polytope $P.$
It could be called 'the index of convex embeddedness of $x$'.
Then, in particular, vertices of $P$ would be of index 1.
However, we refrain from coining a special term.
The common state, for example, in combinatorial optimization, is that
when a polytope is generated by a set of integral points, each of these points is a vertex.
In particular this is true for the travelling salesman polyhedron and other $(0,\, 1)$-polytopes.
Perhaps, this is a reason why the classes of points similar to $C_{\xi}$ were not considered earlier.
\end{RR}

No criterion for $x\in C_{3} (n)$ is known but the computations show that such an $x$ always admits a representation
 $x=\sum \limits _{j=1}^{3}\lambda _{j} y^{j}  ,$ \textit{$y^{j} {\rm \,\vdash\; }n$},
 $\lambda _{j} \ge 0,$ $\sum \limits _{j=1}^{3}\lambda _{j}  =1,$
with all $\lambda _{j} =\frac13.$
The following theorem states that this holds for every integral polytope.
Recall that a polytope is called integral if all its vertices are integer points.

\begin{TT} \label{3points}
If $P\in {\rm R}^{n} $ is an integral polytope and an integer point $x\in P$ is a convex combination
of three integer points in $P$ but is not a convex combination of any two integer points in $P$
then there exist integer points $y^{1} ,\, y^{2} ,\, y^{3} \in P$ such that
\begin{equation} \label{Coef13}
x=\frac13 y_1 +\frac13 y_2 +\frac13 y_3\ .
\end{equation}
\end{TT}

\noindent
\textbf{Proof.} We begin with the general case of an arbitrary integer $k>2$ and an integer $x\in P,$
which is a convex combination of $k$ integer points in $P$
and is not a convex combination of any less than $k$ integer points in $P.$
Then $x$ is a strictly interior point in the $(k-1)$-dimensional simplex $S$ with vertices in these $k$ points.
Assume there is one more integer point $z\in S,$ $z\neq x.$

If $z$ is strictly interior to $S$ then
it divides $S$ to integral simplices $S_1,S_2,\ldots,S_{k}$ with vertices $z$ and any $k-1$ vertices of $S.$
Since $x$ is not a convex combination of any less than $k$ integer points in $S,$
it does not lie in any facet of any $S_j.$
Hence $x$ lies strictly inside one of these $(k-1)$-dimensional simplices, say $x \in S_{1}.$
In the other case, if $z$ lies on the boarder of $S$ let $q$ be the smallest number such that
$z$ is strictly interior to some $q$-dimensional face $F$ of $S.$
Then $z$ divides $F$ to $q+1$ integral simplices $F_1,F_2,\ldots,F_{q+1}.$
This implies that the simplex $S$ can be also divided to $q+1$ integral simplices,
each of whose vertices are vertices of some $F_j$ and the vertices of $S$ not belonging to $F.$
As in the previous case, $x$ lies strictly inside one of these simplices, denote it again by $S_1.$

Applying the same reasoning to $S_{1},$ if it contains an integer point $z_1\neq x,$
we come to a $(k-1)$-dimensional integral simplex $S_{2}\subset S_{1}$ with analogous condition on $x\in S_2.$
After repeating this procedure a finite number of times, we obtain a $(k-1)$-dimensional integral simplex
$T \subset S_{2} \subset S_{1} \subset P$
with $x$ as its single strictly interior integer point satisfying conditions of the theorem
and no integer points on the boarder of $T.$

From here on, we consider that $P$ is the triangle $T$ and $k=3,$ as in the theorem statement.
The rest of the proof can be carried with the help of the Pick's theorem as, for example, in \cite{Rez}.
We will continue using only elementary geometry.
Figure 3 shows the triangle $T$ with vertices $A,B,C$ and the point $x$ denoted by $O.$

      \begin{center}
      \includegraphics[scale=0.4]{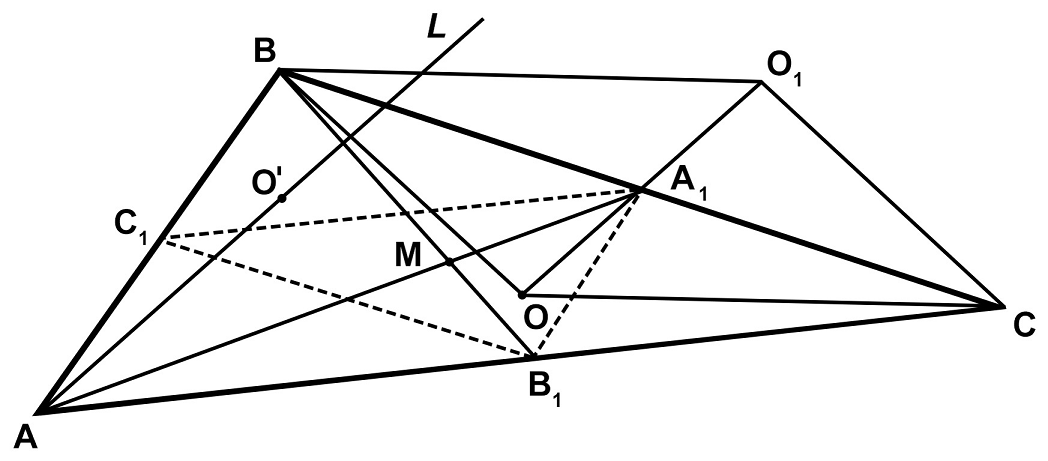}
      \end{center}
      \begin{center}
      {\bf Fig. 3.}  To the proof of Theorem 3.
      \end{center}

We will use the following property of the lattice $H_I$ of integer points in the plain $H$ that contains $T:$
if for some $u_0 \in H_I$ and some $n$-dimensional vector $\bar{c},$
the points $u_0 + \bar{c}$ belong to $H_I$ then
for every $u \in H_I$ the point $u \pm \bar{c}$ belongs to $H_I.$

Let $A_1, B_1, C_1$ be the midpoints of the sides of $T$ and $M$ be the barycenter of $T.$
Assume $O \neq M.$
Then $O$ lies strictly inside the triangle $A_1B_1C_1$
since otherwise, if for example $O\in \triangle A_1B_1C,$
we would have the point $O+\overline{CO} \neq O$ in $\in T\cap H_I.$
(Here $\overline{CO}$ is the vector from $C$ to $O.$)
Hence $O$ lies in one of the triangles $A_1B_1M,$ $A_1C_1M,$ $B_1C_1M$
or on a common side of some two of them.
Let $O\in \triangle  A_1B_1M.$
Draw the parallelogram $BOCO_1$ on the straight line segments $OB$ and $OC.$
By the above property, $O_1\in H_I.$
The diagonal $OO_1$ of the parallelogram passes through $A_1.$
Draw the ray $\mathcal{L}$ parallel to $OO_1$ from $A$ inside the triangle $ABC.$
Since $A_1O$ goes between $A_1M$ and $A_1B_1,$ where $A_1M$ is allowed but $A_1B_1$ is not,
$\mathcal{L}$ goes between $AM$ and $AB$ and can contain $AM$ but not $AB.$
Put the point $O'$ on $\mathcal{L}$ at the distance $|AO'|=|OO_1|$ from $A.$
By the above property, $O'$ is an integer point.
The triangles $ABM$ and $A_1B_1M$ are congruent with the congruence coefficient 2 and $|AO'|=2|A_1O|.$
Hence $O\in \triangle A_1B_1M$ implies $O' \in \triangle ABM.$
Note that $O'$ can lie on $AM$ or $BM.$
In any case $O'$ is in $T$ and integrality of $O'$ implies $O'=O.$

Since $M$ is the single common point in $\triangle A_1B_1M \cap \triangle ABM$
the assumption $O\neq M$ implies $O'\neq O.$
The contradiction proves that $O$ is the barycenter of $T$ and satisfies (\ref{Coef13}).
\hfill $\Box$

Since all integer points in $P_{n}$ are partitions of $n$ \cite{Shl1},
Theorem \ref{3points} implies the following corollary.

\begin{CC}  \label{3Prt}
Every partition $x\in C_{3}(n)$ is the barycenter of some partitions
$y^{1} ,\, y^{2} ,\, y^{3} {\rm \,\vdash\; }n$ (so that the equality (\ref{Coef13}) holds).
\end{CC}

All known partitions $x\in C_{4} (n)$ admit convex representations with coefficients $\tfrac14.$
However the analogue of Theorem \ref{3points} does not hold for such an $x.$
This follows from the results of B.~Reznik \cite{Rez}:
in the case of an integral simplex $P$ with 4 vertices (a 3-dimensional tetrahedron)
there are 7 variants for the values of coefficients
in a convex representation of a single integer point in $P$ via its vertices.
It is interesting that in each variant all denominators are simultaneously equal to one of the numbers 4, 5, 7, 11, 13, 17, 19.
Nothing is known about the coefficients in convex representations via 5 integer points.

\section{Evenness of $n$ and Metropolis partitions}

Let us return to Figure 1 that presents the graph of the function $v(n)$ for $n\leq 100.$
One immediately sees that the value of $v(n)$ depends on the evenness of $n$:
\begin{equation} \label{GrindEQ__2_}
v(2r-1)>v(2r)
\end{equation}
except small $r.$ So, we can refer to the $v(n)$ graph as consisting of two subgraphs:
for odd and even $n$'s, the latter lying below the former.
This radically differs from the monotone increasing of $p(n),$ the number of partitions of $n.$

Upon careful examination of Figure 1 we suspected that some points $(n,v(n))$ with $n$ odd
are disposed slightly higher than the main line.
It turned out that they correspond to prime $n$'s.
Comparison of their heights $v(n)$ with the half-sums $\tfrac12(v(n-2)+v(n+2))$
confirmed this observation for all prime $n\ge 43$ except $n=61.$
The observed tooth-shaped form of the $v(n)$ graph
and special role of prime numbers raised the question of what multiplicative property of $n$ affects the value of $v(n).$

We know from the computation that for every $n,$ the majority of partitions that are not vertices belong to $C_{2} (n).$
By Theorem 2, these partitions have two collections of parts with equal sums.
In particular, for even $n=2r,$ $C_{2}(n)$ contains partitions of the form
\begin{equation} \label{Metro}
 [\text{partition\_1 of~} r] + [\text{partition\_2 of~} r],
\end{equation}
where
\begin{equation} \label{Metro_1}
 \text{partition\_1} \neq \text{partition\_2}\ .
\end{equation}

Denote the number of partitions (\ref{Metro}), disregarding condition (\ref{Metro_1}), by $m_{2} (2r).$
It is not hard to see that
\begin{equation} \label{MetroBound}
 m_{2}(2r)=\frac12\left(\, p(r)^{2} +p(r)\right) - [\text{number of duplicates in (\ref{Metro})}],
\end{equation}
but it is far from clear how to count the duplicates.
Note that if a partition of the form (\ref{Metro}) satisfies (\ref{Metro_1}) it can be a vertex.
The partition $(0,2,2,0^7)\vdash 10$ is an example.
The number of such vertices is less than $p(r),$ which is a rough estimate.
Thus, when we are interested in the asymptotics of $C_{2}(n)$ and $v(n),$
we can ignore vertices of the form (\ref{Metro}).
Note that these vertices are shown in Figure 1 by the small rectangle in the ${\rm Vert}\, P_{n}$ area.
The following conjecture may explain inequality (\ref{GrindEQ__2_}) and the tooth-shaped form of the $v(n)$ graph.

\begin{Conj} \label{Conj2}
For $n$ even, $m_{2} (n)$ is large relative to $v(n).$
\end{Conj}

After having searched in the OEIS by the sequence of the first values of $m_{2}(n),$
we have got to the sequence A002219 and the work of Metropolis and Stein \cite{Metro},
where the authors had counted partitions of $n$ that can be obtained by joining $r,$ $r$ divides $n,$
not necessarily different partitions of $\frac{n}{r}$
(for convenience, we slightly changed the original notation in \cite{Metro}).
We call these partitions Metropolis $r$-partitions.
For $n$ multiple of $r,$ denote the set of Metropolis $r$-partitions of $n$ by $M_{r} (n)$ and set $m_{r} (n):=|M_{r} (n)|.$
Note that Metropolis 2-partitions coincide with partitions (\ref{Metro}).
The main result of \cite{Metro} is the formula for $m_{r}(n)$
in the form of a finite series of binomial coefficients multiplied by certain integer coefficients,
which depend only on $r.$
For $m_{2}(n)$ this formula reads
\begin{equation} \label{m2n}
m_{2}(n)={{g+2}\choose{2}}+(g+2)c_{1} + c_{2},
\quad g=\bigg\lfloor \frac{\tfrac{n}{2}+1}{2} \bigg\rfloor, \quad
\dfrac{n}{2}>5,
\end{equation}
where $c_{1} $ and $c_{2} $ 'must be determined by direct calculation' \cite{Metro}.
The sequence A002219 contains the values of $m_{2} (n)$ for even $n\le 178.$
Using (\ref{m2n}), we obtain an upper bound $b(n)$ for the number of vertices of $P_n.$

\begin{TT} \label{Bound}
\begin{equation}\label{b_n}
v(n)\leq b(n) :=
\begin{cases}
p(n)-m_{2}(n),   & n \text{~even},\\
p(n)-m_{2}(n-1), & n \text{~odd},
\end{cases}
\end{equation}

\noindent
where values of $m_{2}(\cdot)$ are calculated with the use of (\ref{m2n}).
\end{TT}

\noindent
\textbf{Proof.} The proof follows from the inclusion $M_{2}(n) \subset C_{2}(n),$
if we ignore the small number of vertices belonging to $M_{2}(n),$
and the fact that adding the part 1 to every partition in $M_{2}(n-1),$ $n$ odd,
results in a partition in $C_{2}(n).$ \hfill $\Box$

Disregarding the duplicates in (\ref{Metro}), one can obtain from (\ref{MetroBound}) an upper bound on $m_2(n).$
However, Metropolis and Stein pointed that, for large $n,$ much better is the bound
$p(n,\frac{n}{2}),$ which is the number of partitions of $n$ with no part greater than $\frac{n}{2}.$
It is not hard to shaw that $p(n,\frac{n}{2})$ is asymptotically equal to $p(n).$
An anonymous author under the nickname 'joriki' presented the following proof of this fact in Stackexchange \cite{joriki}.
Every partition of $n$ has at most one part $m$ larger than $\frac{n}{2},$ and the remaining parts form a partition of $n-m.$
Thus
$$
p\left(n,\frac{n}{2}\right) = p(n) - \sum\limits_{i=0}^{\frac{n}{2}-1}p(i).
$$
For large $n,$ the terms in the sum are exponentially smaller than  $p(n),$ so asymptotically
$$
p\left(n, \frac{n}{2}\right) \sim p(n).
$$

One can draw the graph of $m_2(n)$ over the known values and see that the equivalence
\begin{equation}\label{m2equiv}
    m_2(n) \sim p(n)
\end{equation}
is also very likely to be true.

Let us turn to the relations between $m_{2}(n),$ $p(n),$ and $v(n),$ whose values are known for $n\le 100.$
Table 1 presents some data for a few values of $n.$
We see that the ratios $v(n)/p(n)$ and $v(n)/m_{2} (n)$ are small and very close,
and the ratio $\left(p(n)-m_{2} (n)\right)/p(n)$ rapidly decreasses when $n$ grows,
though the difference $p(n)-m_{2} (n)$ also increasses.
This corroborates (\ref{m2equiv}) and gives an additional argument in favor of Conjecture \ref{Conj2}.

\footnotesize
\begin{table}[h]
\centering
\begin{tabular}{|*{8}{r|}} \hline
\small$n$ & \small$p(n)$ & \small$v(n)$ & \small$m_{2} (n)$ & \small$\frac{v(n)}{p(n)}$ & \small$\frac{v(n)}{m_{2}(n)}$ &
\small$p(n)-m_{2}(n)$ & \small$\frac{p(n)-m_{2}(n)}{p(n)}$ \\ \hline
60 & 966467 & 5148 & 924522 & 0.005327 & 0.005568 & 41945 & 0.0434 \\ \hline
78 & 12132164 & 17089 & 11850304 & 0.001409 & 0.001442 & 281860 & 0.0232 \\ \hline
100 & 190569292 & 59294 & 188735609 & 0.000311 & 0.000314 & 1833683 & 0.0096 \\ \hline
\end{tabular}
\caption{Relations between $p(n),$ $v(n),$ and $m_{2}(n)$} \label {tab : metka}
\end{table}
\normalsize

The expression (\ref{b_n}) for $b(n),$ the upper bound on the number of vertices of $P_n,$ may help to
clarify, though not prove, the cause of the tooth-shaped form of the graph of $v(n)$ under Conjecture 2.
For $n$ odd, it yields
$$
b(n) - \frac12\Big(b(n-1) + b(n+1)\Big) = \Big(p(n) - \frac12(p(n-1) + p(n+1))\Big) + \frac12\Big(m(n+1) - m(n-1)\Big),
$$
where the first term is asymptotically zero and the second term is positive.
This means that $b(n)$ has a peak at every large odd $n$
and the graph of $b(n)$ is of the tooth-shaped form, similar to that in Figure 1 for $v(n).$

Let us consider two examples to see what happens when $n$ is even. For $n=78,$ we have
$b(78) = p(78)-m_{2}(78)=281\,860$, while
$b(77) = p(77)-m_{2}(76)=1\,549\,719.$
So, the bound for $v(78)$ is less than $0.19 \cdot b(77).$
In the same way we have $b(100) < 0.09 \cdot b(99)\,$!
Hence it is more than likely that, for $n$ even, $v(n)$ is not only less than $\frac12(v(n-1) + v(n+1))$ but
$v(n)<v(n-1).$
Thus, Conjecture \ref{Conj2} and the asymptotic equivalence (\ref{m2equiv}),
as its stronger form observed from the numerical data, reasonably justify the inequality (\ref{GrindEQ__2_}) and
the gap between the values of $v(n)$ for even and odd $n$'s.

The following theorem provides a supplemental indication of the importance
of Metropolis 2-partitions for recognizing vertices of $P_n.$
Call a partition $x{\rm \,\vdash\; }n$ an extension of a partition $y{\rm \,\vdash\; }m$, $m<n$,
if every part of $y$ is a part of $x.$

\begin{TT} \label{Metro2}
For every $n,$ every partition $x\in C_{2} (n)$ is either a Metropolis 2-partition
or an extension of some Metropolis 2-partition $y{\rm \,\vdash\; }m,$ $m<n.$
\end{TT}

\noindent
\textbf{Proof.}
Consider arbitrary $n$ and $x\in C_{2}(n),$ $x\notin M_{2} (n)$ if $n$ is even.
By Theorem \ref{Criterion}, there exist two collections of parts of $x$ with the same sum.
Let $s$ be the minimal value of such a sum. Clearly, $s\leq\frac{n}{2}.$
The corresponding collections are disjoint and their union is a Metropolis 2-partition $y$ of some $m=2s\leq n.$
Hence $x=y$ if $m=n$ or $x$ is an extension of $y$ if $m<n.$
\hfill $\Box$

\section{$n$'s multiple of 3 and knapsack partitions}

R. Ehrenborg and M. A. Readdy \cite{EhRe07} called a partition $x$ a knapsack partition if for every integer,
there is utmost one way to represent it as a sum of some parts of $x.$
Denote the set of knapsack partitions of $n$ by $K(n)$ and set $k(n):=|K(n)|.$
Theorem \ref{Criterion} implies relations
\[
K(n)=P(n)\backslash C_{2} (n),
\]
\[
C_{\xi } (n)\subset K(n),\, \, \, \, \, \, \xi >2,
\]
\[
{\rm Vert}\, P_{n} \subset K(n).
\]
The smallness of $|{\rm Vert}P_{n}\cap M_2(n)|$ implies that for large $n,$
\[
v(n)<k(n)<p(n)-m_{2}(n).
\]

Note that $k(n)$ is a much better upper bound on $v(n)$ than $b(n)$ in (\ref{b_n})
but no formula for $k(n)$ is known.

Ehrenborg and Readdy computed the values $k(n)$ for $n\le 50$ and exhibited them in the OEIS, sequence A108917.
We extended this sequence till $n=165$ as a by-product of our computation of vertices of $P_{n}.$
Table 2 enhances Table~1 by the $k(n)$ values.
Consider its first three rows with even $n.$
Looking at the columns ${v(n)\mathord{\left/ {\vphantom {v(n) (p(n)-m_{2} (n))}} \right. \kern-\nulldelimiterspace} (p(n)-m_{2} (n))}$
and $k(n)-v(n)$ and comparing the columns $k(n)$ and $p(n)-m_{2} (n),$ we see that
many partitions of $n$ that are neither vertices of $P_{n}$ nor Metropolis 2-partitions
are convex combinations of 2, 3, or more partitions of $n.$

\scriptsize
\begin{table}[h]
\centering
\small
\begin{tabular}{|*{9}{r|}} \hline
\scriptsize$n$ & \scriptsize$p(n)$ & \scriptsize$v(n)$ & \scriptsize$m_{2}(n)$ &
\scriptsize$p(n)-m_{2}(n)$ & \footnotesize$\frac{v(n)}{p(n)-m_{2}(n)}$ &
\scriptsize$k(n)$ & \footnotesize$\frac{v(n)}{k(n)}$ & \scriptsize$k(n)-v(n)$ \\ \hline
60 & 966\,467 & 5\,148 & 924\,522 & 41\,945 & 0.12 & 5 341 & 0.964 & 193 \\ \hline
78 & 12\,132\,164 & 17\,089 & 11\,850\,304 & 281\,860 & 0.06 & 17\,871 & 0.956 & 782 \\ \hline
100 & 190\,569\,292 & 59\,294 & 188\,735\,609 & 1\,833\,683 & 0.03 & 61\,692 & 0.967 & 2\,398 \\ \hline
77 & 10\,619\,863 & 21\,393 &  &  &  & 22\,128 & 0.967 & 735 \\ \hline
\end{tabular}
\normalsize
\caption{Relations between $p(n)$, $v(n)$, $m_{2}(n),$ and $k(n)$} \label {tab : metka_1}
\end{table}
\normalsize

One can check that the graph of $k(n),$ like that of $v(n),$ disintegrates into two graphs, for $n$ odd and $n$ even.
However we see that the ratio $v(n)/k(n)$ does not increase monotonically and
is approximately the same for odd $n=77$ and even $n=100,$ which are rather far from each other.
Figure 4 presents the graph of $v(n)/k(n).$
We obviously see that the $n$'s multiple of 3 are the local minima of $v(n)/k(n).$
This means that such $n$'s have more partitions that are not vertices of $P_n$
and do not belong to $C_{2} (n)$ than the $n$'s not multiple of 3.


\begin{center}
\noindent \includegraphics*[width=6in, height=1.6in, keepaspectratio=false]{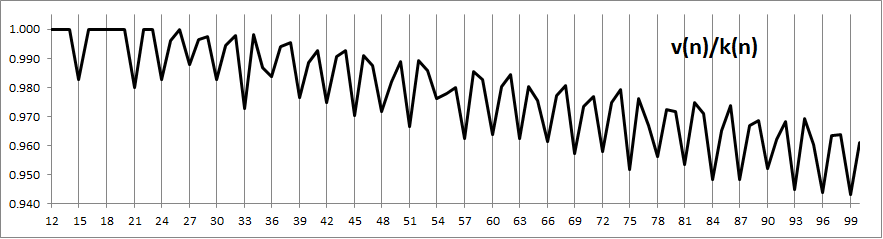}
\end{center}
      \begin{center}
      {\bf Fig. 4.} Ratio $v(n)/k(n)$ of the number of vertices
                    to the number of knapsack partitions.
      \end{center}


The following conjecture naturally explains this phenomenon.

\begin{Conj} \label{Conj3}
The majority of partitions in $K(n)$ that are not vertices of $P_{n} $
are convex combinations of three partitions of $n.$
\end{Conj}

Conjecture \ref{Conj3} is consistent with our computation experience.
We know that for $n$ multiple of 3, most partitions in $C_{3}(n),$ that are not
extensions of some partitions in $C_{3}(q),$ $q<n,$ $q$ multiple of 3, with the additional part $n-q,$
have a part $\frac{n}{3}$ and
one of the partitions involved in the convex combination has three parts $\frac{n}{3}.$
For example, the non-vertex $x=(1^{3} ,9,17,22){\rm \,\vdash\; }51$ has a part
$17=\frac{51}{3},$ and its convex representation is
$\frac{1}{3}(1^{7} ,22^{2} ) + \frac{1}{3}(1^{2} ,9^{3} ,22) + \frac{1}{3}(17^{3}).$

As for the tendency of $v(n)/k(n)$ to decrease, we see its explanation
in the increase of the number of partitions in $C_{\xi } (n)$, $\xi >3,$ with the growth of $n.$

\section{Stratification of the numbers of vertices}

To reveal the discovered dependence of the number of vertices of the polytope $P_n$ on multiplicative properties of $n$
and examine it in more details, we consider the classes of integers
$$
N_{k} :=\left\{n\ |\ n=kp,\, \, p~ \text{prime}, k\leq p \right\},\, \, \, \, k=1,2,3,...
$$
and the corresponding numbers of vertices
$$
v_{k} (n):=v(n),\, \, \, n\in N_{k}.
$$

Figure 5 demonstrates the graphs of the functions $v_{k} (n)$ for $k=1,\, 7,\, 5,\, 3,\, 2,\, 4,\, 6$
in order from top to bottom.
They are generated with the use of the FindFit method of Wolfram Mathematica.
We approximated the known values of $v_{k}(n)$ by the functions of the form $Ae^{ B\sqrt{n}}$ with parameters $A$ and $B.$
The segment $n\in [60, 70]$ is chosen to split the graphs $v_{k}(n)$ visually.
It also lies in the most interesting part of the segment $[1,100],$
where we can expect our approximations to reveal a reliable picture of what happens.
We do not include the graphs of $v_{k}(n),$ $k>7,$ in Figure 5
because they are little informative.
For these $k,$ there are too few prime numbers $p\geq k$ such that $kp \in [1, 100].$
We also disregard the condition $k\leq p$ for $k=7$ when drawing the graph of $v_{7}(n),$
since it holds only for $n=49$ and $n=77.$

\begin{center}
\noindent \includegraphics*[width=6.45in, height=2.50in,
keepaspectratio=false]{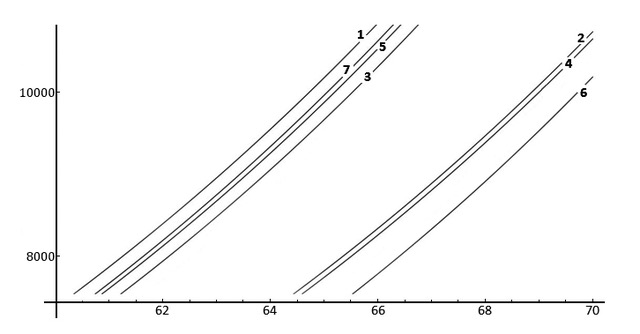}
\noindent {\bf Fig. 5.} Stratification of the number of vertices of the polytope $P_{n}$ function \\
to the functions $v_k(n),$ $k=1,2,\ldots,7.$
\end{center}


We see that the graph of $v(n)$ is neither a single line nor a conjunction of two lines,
for odd and even $n,$ as in Figure 1.
It is stratified into layers corresponding to the classes $N_{k}$ and resembles a layered cake.
Its layers are of the same shape but are disposed at different levels.
The topmost line corresponds to $N_{1},$ the class of primes.
The graph of $v_{3}(n)$ goes below it.
Between them, one below another, are disposed the graphs of $v_{7}(n)$ and $v_{5}(n),$
while for $k$ even, the graphs of $v_{k}(n)$ go below $v_{3}(n).$

The levels of the graphs of $v_{1}(n),$ $v_{2}(n),$ $v_{3}(n)$ agree with Conjectures 2 and 3.
The intermediate position of $v_{7}(n)$ and $v_{5}(n)$ and Conjecture 1
move us to suggest a more general conjecture that for prime $k$ dividing $n,$
the determining influence on the level of $v_{k}(n)$ is exerted by
the number of partitions of $n\in N_{k}$ that belong to $C_{k} (n).$
If $k_{1}, k_{2}$ are two primes, $k_{1} >k_{2},$ then,
in accordance with Conjecture 1, for large and sufficiently close to each other
$n_{1} \in N_{k_{1} } $ and $n_{2} \in N_{k_{2} },$ the inequality
$|C_{k_{1}}(n_{1})| < |C_{k_{2}}(n_{2})|$ holds
and therefore the $v_{k_1}(n)$ graph is disposed above the $v_{k_2}(n)$ graph.

The case of $k,$ a composite divisor of $n,$ can be explained
using the graphs of $v_{6}(n)$ and $v_{4}(n).$
$v_{6}(n)$ is disposed below $v_{2}(n)$ and $v_{3}(n)$ because the level of
$v_{6}(n)$ is affected by partitions in $C_{2}(n)$ and partitions in $C_{3}(n).$
Similarly, $v_{4}(n)$ goes between $v_{2}(n)$ and $v_{6}(n)$ because
$4$ is an additional (to $2$) divisor of $n$ and $|C_{4}(n)|<|C_{3}(n)|.$

We summarize the above in the final conjecture.

\begin{Conj}\label{ConjFinal}
The number of vertices of $P_n$ depends on factorization of $n.$
The graph of $v(n)$ is stratified into layers $v_k(n).$
This stratification is based on partitioning of integer numbers to the classes $N_k.$
Let $n\in N_k,$ $1\leq k\leq p,$ $p$ prime, and let $k_{1} ,k_{2} ,k_{3} , \ldots$ be the divisors of $k$ sorted in ascending order.
Then the major influence on the level of the graph of $v_k(n),$ on which the value of $v(n)$ lies, is rendered by the divisor $k_{1}.$
Every successive divisor makes its additional contribution to lowering the level of $v_k(n).$
\end{Conj}

Generalizing, we might say that the value of $v(n)$ is determined by
the proximity of $n$ to its greatest prime divisor,
which is defined by the lexicographic order on the set of increasing sequences of divisors of $n$.
For example, $38=2\cdot 19$ would be 'more prime` than $39=3\cdot 13,$
hence the layer $v_2(n),$ that contains $v(38),$ is disposed lower than the layer $v_3(n)$ containing $v(39).$
The same would hold for $78=2\cdot 3\cdot 13$ and $70=2\cdot 5\cdot 7.$
If we extend this speculation, we might come to a fractal structure of the graph of $v(n).$
For example, the graph of $v_5(n)$ together with $v_{10}(n), v_{15}(n), v_{20}(n), \ldots$
may have a structure similar to that of $v(n).$
However, it is too early to foresee so far ahead --- more numerical data of $v(n)$ is needed.
Then the Conjecture 4 might be further detalized.

\section{Remark on the Gomory's corner polyhedron}

Let $G$ be a finite Abelian group, $G^+$ be the set of its nonzero elements, and $g_0\in G.$
The master corner polyhedron $P(G,g_0)$ was defined by R. E. Gomory \cite{Gom} as
the convex hull of solutions
$$t=(t(g)\ ;\ g\in G^+)\in \mathbb{R}^{|G^+|},\ ~t(g)~ \text{integer}, ~t(g)\geq 0,$$
to the equation
\begin{equation}\label{groupEq}
 \sum_{g\in G^+}t(g)g = g_0.
\end{equation}
For $G_{n+1} := \mathbb{Z}/(n+1)\mathbb{Z},$ the cyclic group of order $n+1,$ and $g_0=n,$
the equation (\ref{groupEq}) reads
\begin{equation*}
 t_1 + 2t_2+ \ldots +nt_n \equiv n\bmod(n+1),
\end{equation*}
which differs from (\ref{eqPn}) only in that the addition here is modulo $n+1.$
That is why we call $P(G,g_0)$ the elder brother of $P_n.$

Our experience in studying both polyhedra shows that the vertex structure of $P_n$
is more transparent and easy for understanding than that of the $P(G,g_0),$ even in the case of the cyclic group.
In our opinion, this is because the standard addition on the segment of integers $[1, n],$
albeit defined only partially, is much easier to comprehend than the group addition.
Most results on vertices of $P_n$ were successfully transferred to vertices of the master corner polyhedron \cite{EJOR}.

Statistics on vertices of $P(G,g_0)$ is unbelievably poor.
For many years, all that we knew about their numbers could be found in the R. E. Gomory's seminal paper \cite{Gom}.
The researchers concentrated their efforts on studying facets of $P(G,g_0)$
since they induce the most efficient cuts for the integer linear programs.
In contrast, vertices --- though they are no less important for understanding the structure of $P(G,g_0)$ ---
fell out of research.
In~\cite{Gom} Gomory computed vertices of $P(G,g_0)$ for all groups $G$ of the order up to 11 and all $g_0 \in G.$%
\footnote{One extra point $t\in P(G_{11},10),$ with $t(5)=1,$ $t(9)=3$ and all other $t(i)=0,$
indicated in \cite{Gom} as a vertex was excluded in \cite{EJOR}.}
The numbers of these vertices for the the case of corner polyhedra $P(G_{n+1},n),$ $n=1,2,\ldots,10,$
constitute the sequence A300795 in the OEIS.
Recently, D. Yang extended this data till $n=21$ \cite{Yang}.
Figure 6 exhibits the graph of the final sequence.

\begin{center}
\noindent \includegraphics*[width=6.2in, height=2.50in,
keepaspectratio=false]{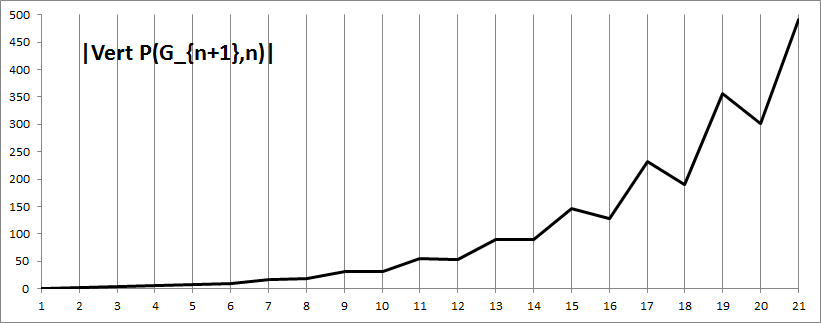}
\noindent {\bf Fig. 6.} Graph of the number of vertices function
for the master corner polyhedron \\ $P(G_{n+1},n)$ on the cyclic group $G_{n+1},$ $1\leq n\leq 21.$
\end{center}

We perceive this picture as a forerunner of a graph similar to that depicted in Figure 1.
The tooth-shaped form of the $|{\rm Vert}\, P(G_{n+1},n)|$ graph is obvious even in this initial part.
Some of the above features of $v(n)$ may also become visible
when the sequence of numbers of vertices of $P(G_{n+1},n)$ will be extended.

\section{Concluding remarks}

In order to study the number of vertices $v(n)$ of the polytope $P_n$ of integer partitions of $n$
we investigated the structure of the set of partitions that are not vertices.
We divided this set to disjoint subsets $C_{\xi}(n)$ according to the minimum number $\xi$ of partitions
needed to represent a partition as their convex combination.
Using the available numerical data, we demonstrated that
$M_2(n),$ the set of Metropolis 2-partitions of $n,$
constitutes a larger part of partitions that are not vertices of $P_n.$
As a consequence, vertices of $P_n$ form a small subset of partitions of $n.$
We proved that an integer point in an arbitrary integral polytope $P,$
which belongs to the subset of integer points in $P$ analogous to $C_3(n),$
admits a convex representation via three integer points with all coefficients equal to $\frac13.$

Thorough analysis of the computed values of $v(n)$ revealed intriguing properties of this function.
Comparing this data with available numbers of knapsack and Metropolis 2-partitions
moved us to suggest several conjectures that explain observed peculiarities.
The main conjecture claims that $v(n)$ depends on factorization of $n.$
We presented visual but convincing arguments in its favor.
We showed that the graph of $v(n)$ is stratified into layers, the subgraphs corresponding to the classes of integers
that are determined by factorization of $n.$
The upper layer corresponds to prime numbers and the others correspond to collections of small divisors of $n.$
Every prime divisor makes its own contribution to lowering the level of the layer.
The smaller the divisor the more significant its effect.

We provided an argument in favor of a similar dependence
for the number of vertices of the master corner polyhedron on the cyclic group.
Though the data in our disposal is still rather limited, we believe that this argument deserves further examination
in view of the closeness of the master corner polyhedron and $P_n$ revealed in \cite{EJOR}.

This work draws forth new questions.
Formal proof and detailed study of the dependence of the number of vertices of $P_{n}$
on factorization of $n$ remain open problems for the future research.
Further computation of $v(n)$ would be of great help.
One of the most important problems is to find a combinatorial criterion for vertices of $P_n.$
More specific problems are concerned with the nature of partitions in $C_{\xi}(n).$
Counting knapsack partitions does not look unworkable.
This problem looks easier than enumerating vertices.
Its solution will provide a rather good estimate for $v(n).$
We also hope that this work will give an impetus to further study of vertices of the corner polyhedron.

\section{Acknowledgements}

We are deeply obliged to the late Professor V.~B. Priezzhev of the Joint Institute for Nuclear Research (JINR), Dubna, Russia,
who has recently passed away, for useful discussions.
His friendly and encouraging support will be greatly missed.
This work could not have been done without computation of vertices performed by A.~S. Vroublevski.
We thank D. Yang, a student at UC Davis, for computing vertices of the master corner polyhedron.
We also thank E. S. Zabelova for her assistance in applying Wolfram Mathematica.

\end{document}